 \documentclass[12pt,leqno]{article} 
\usepackage{euscript} \def\mathcal{\EuScript}
 \def\DATE{December  9, 2004} 
 \usepackage{amssymb} 
 \newtheorem{theorem}{Theorem} 
 \newtheorem{definition}[theorem]{Definition}

 \newtheorem{example}[theorem]{Example}

 \newtheorem{proposition}[theorem]{Proposition} 
 \newtheorem{remark}[theorem]{Remark}

 \catcode`\@=11 

 \def\ps@myheadings{\let\@mkboth\@gobbletwo 
 \def\@oddhead{\ifnum\count0=1 \hfill\else 
 \rightmark \hfil \rm\thepage\fi}%
 \def\@oddfoot{\ifnum\count0=1 \hfill \rm 1 \hfill \else 
 \hfill\fi} 
 \def\@evenhead%
 {\rm\leftmark\hfil\rm\thepage}%
 \def\@evenfoot{}\def\sectionmark##1{} 
 \def\subsectionmark##1{}} 

 \def\@begintheorem#1#2{\it \trivlist \item[\hskip 
  \labelsep{\bf #1\ #2.}]} 
 \def\@opargbegintheorem#1#2#3{\it \trivlist\item[\hskip%
  \labelsep{\bf #1\ #2.\ (#3)}]} 
 \def\@endtheorem{\endtrivlist} 

 \def\@listI{\leftmargin\leftmargini \parsep 1pt plus 2.5pt 
  minus 1pt\topsep 5pt plus 2pt minus 3pt%
  \itemsep 0pt plus 2.5pt minus 1pt} 
 \let\@listi\@listI 
 \@listi 

 \def\@sect#1#2#3#4#5#6[#7]#8{\ifnum #2>\c@secnumdepth%
  \def \@svsec {}\else \refstepcounter {#1}\edef \@svsec%
  {\csname the#1\endcsname. \hskip .1em }\fi \@tempskipa%
  #5\relax \ifdim \@tempskipa >\z@ \begingroup #6\relax%
  \@hangfrom {\hskip #3\relax \@svsec }{\interlinepenalty%
  \@M #8\par }\endgroup \csname #1mark\endcsname {#7}%
  \addcontentsline {toc}{#1}{\ifnum #2>\c@secnumdepth%
  \else \protect \numberline {\csname the#1\endcsname. }%
  \fi #7}\else \def \@svsechd {#6\hskip #3\@svsec #8%
  \csname #1mark\endcsname {#7}\addcontentsline {toc}{#1}%
  {\ifnum #2>\c@secnumdepth \else \protect \numberline%
  {\csname the#1\endcsname. }\fi #7}}\fi \@xsect {#5}} 

 \def\section{\@startsection {section}{1}{\z@ }%
  {-3.5ex plus -1ex minus -.2ex}{2.3ex plus .2ex}{\bf }} 

 \def\thebibliography#1{%
  \section *{References\@mkboth {REFERENCES}{REFERENCES}}%
  \list {[\arabic {enumi}]}{\settowidth \labelwidth {[#1]}%
  \leftmargin \labelwidth \advance \leftmargin \labelsep %
  \usecounter {enumi}} \def \newblock %
  {\hskip .11em plus .33em minus -.07em} \sloppy \clubpenalty 4000%
  \widowpenalty 4000 \sfcode`\.=1000\relax} 

 \def\@maketitle{%
  \newpage \null \vskip 2em 
  \begin{center} 
 {\Large\bf \@title \par } 
  \vskip 1.5em 
  {\large \lineskip .5em 
  \begin {tabular}[t]{c}\@author 
  \end{tabular}\par} 
  \end{center} 
   \vskip .8em}

 \def\abstract{%
 \if@twocolumn \section *{Abstract} 
  \else \small\quotation\noindent{\bf Abstract.}\fi} 

 \catcode`\@=13 

\def\Lin{{\rm Lin}} \def\Yhoch{{Y_{\rm Hoch}}} \def\Ycohoch{{Y_{\rm coHoch}}}
\def\otexp#1#2{{#1^{\otimes #2}}}
\def\Lie{{\mathcal{L}{\it Lie}}}
\def\pre-Lie{{\mbox {\it pre-\Lie}}}
\def\Ass{{\mathcal{A}{\it ss}}}
\def\Gass#1{{\mbox {$\mathcal{G}_#1$\it-ass}}}
\def\Vinb{\mathcal{V}{\it inb}}
 \def\Poiss{\mathcal{P}{\it oiss}}
 \def\Com{\mbox{{$\mathcal C$}\hskip -.3mm {\it om}}} 
 \def\BOX{\raisebox{1.2mm}{\hskip .5mm $\fbox{\hphantom{\hglue .01mm}}$}} 
 \def\Vect{{\tt Vect}} 
 \def\Palg{\mbox{{$\mathcal P$}-{\tt alg}}} 
 \def\calP{{\mathcal P}} 
  \def\KK{{\mathbb {K}}}
 \def\id{1\!\!1} 
 \def\sgn{{\rm sgn}} 
 \def\ot{\otimes} 
 \def\End{\hbox{${\mathcal E}\hskip -.1em {\it nd}$}}
 \def\Lie{\hbox{{$\mathcal L$}{\it ie\/}}}
 \def\qed{\hspace*{\fill}
    \mbox{\hphantom{mm}\rule{0.25cm}{0.25cm}}\\}

 \def\jedna#1#2{
    \unitlength .5mm
    \begin{picture}(15,9)(6,11)
    \thicklines
    \put(13,17){\line(0,1){4}}
    \qbezier(12,15)(10,13)(8,11)
    \qbezier(14,15)(16,13)(18,11)
    \put(13,15.25){\makebox(0,0)[cc]{\Large $\bullet$}}
    \put(8,9){\makebox(0,0)[tc]{\rm \scriptsize #1}}
    \put(18,9){\makebox(0,0)[tc]{\rm \scriptsize #2}}
\end{picture}
}

 \def\dva#1#2{
    \unitlength .5mm
    \begin{picture}(15,9)(6,11)
    \thicklines
    \put(13,18){\line(0,1){3}}
    \qbezier(11,14)(10,13)(8,11)
    \qbezier(15,14)(16,13)(18,11)
    \put(13,15.25){\makebox(0,0)[cc]{\LARGE $\circ$}}
    \put(8,9){\makebox(0,0)[tc]{\rm \scriptsize #1}}
    \put(18,9){\makebox(0,0)[tc]{\rm \scriptsize #2}}
\end{picture}
}

 \def\tri#1#2{
    \unitlength .5mm
    \begin{picture}(15,9)(6,11)
    \thicklines
    \put(13,18){\line(0,1){3}}
    \qbezier(11,14)(10,13)(8,11)
    \qbezier(15,14)(16,13)(18,11)
    \put(11,18){\line(1,0){4}}
    \put(11,18){\line(0,-1){4}}
    \put(15,14){\line(-1,0){4}}
    \put(15,14){\line(0,1){4}}
     \put(8,9){\makebox(0,0)[tc]{\rm \scriptsize #1}}
    \put(18,9){\makebox(0,0)[tc]{\rm \scriptsize #2}}
\end{picture}
}

\def\dvadva#1#2#3{
    \unitlength .5mm
    \begin{picture}(15,9)(6,11)
    \thicklines
 \put(13,18){\line(0,1){3}}
    \qbezier(11,14)(10,13)(2,6)
    \qbezier(15,14)(15.5,13.5)(16,13)
    \put(13,15.25){\makebox(0,0)[cc]{\LARGE $\circ$}}
    \put(2,4){\makebox(0,0)[tc]{\rm \scriptsize #1}}
\put(5,-5){
    \qbezier(11,14)(10,13)(8,11)
    \qbezier(15,14)(16,13)(18,11)
    \put(13,15.25){\makebox(0,0)[cc]{\LARGE $\circ$}}
    \put(8,9){\makebox(0,0)[tc]{\rm \scriptsize #2}}
    \put(18,9){\makebox(0,0)[tc]{\rm \scriptsize #3}}
}
\end{picture}
}

\def\AAdva#1#2#3{
    \unitlength .5mm
    \begin{picture}(15,9)(6,11)
    \thicklines
 \put(13,18){\line(0,1){3}}
    \qbezier(11,14)(10,13)(2,6)
    \qbezier(15,14)(15.5,13.5)(16,13)
    \put(13,15.25){\makebox(0,0)[cc]{\Large $\bullet$}}
    \put(2,4){\makebox(0,0)[tc]{\rm \scriptsize #1}}
\put(5,-5){
    \qbezier(11,14)(10,13)(8,11)
    \qbezier(15,14)(16,13)(18,11)
    \put(13,15.25){\makebox(0,0)[cc]{\LARGE $\circ$}}
    \put(8,9){\makebox(0,0)[tc]{\rm \scriptsize #2}}
    \put(18,9){\makebox(0,0)[tc]{\rm \scriptsize #3}}
}
\end{picture}
}

\def\dvaAA#1#2#3{
    \unitlength .5mm
    \begin{picture}(15,9)(6,11)
    \thicklines
 \put(13,18){\line(0,1){3}}
    \qbezier(11,14)(10,13)(2,6)
    \qbezier(15,14)(15.5,13.5)(16,13)
    \put(13,15.25){\makebox(0,0)[cc]{\LARGE $\circ$}}
    \put(2,4){\makebox(0,0)[tc]{\rm \scriptsize #1}}
\put(5,-5){
    \qbezier(11,14)(10,13)(8,11)
    \qbezier(15,14)(16,13)(18,11)
    \put(13,15.25){\makebox(0,0)[cc]{\Large $\bullet$}}
    \put(8,9){\makebox(0,0)[tc]{\rm \scriptsize #2}}
    \put(18,9){\makebox(0,0)[tc]{\rm \scriptsize #3}}
}
\end{picture}
}

\def\AAAA#1#2#3{
    \unitlength .5mm
    \begin{picture}(15,9)(6,11)
    \thicklines
 \put(13,18){\line(0,1){3}}
    \qbezier(11,14)(10,13)(2,6)
    \qbezier(15,14)(15.5,13.5)(17,12)
    \put(13,15.25){\makebox(0,0)[cc]{\Large $\bullet$}}
    \put(2,4){\makebox(0,0)[tc]{\rm \scriptsize #1}}
\put(5,-5){
    \qbezier(11,14)(10,13)(8,11)
    \qbezier(15,14)(16,13)(18,11)
    \put(13,15.25){\makebox(0,0)[cc]{\Large $\bullet$}}
    \put(8,9){\makebox(0,0)[tc]{\rm \scriptsize #2}}
    \put(18,9){\makebox(0,0)[tc]{\rm \scriptsize #3}}
}
\end{picture}
}

\begin{document} 
 \bibliographystyle{plain} 

 \title{Algebras with one operation including Poisson and other 
 Lie-admissible algebras} 

 \author{Martin Markl\thanks{Supported by the grant  GA \v CR 
                   201/02/1390.} 
 \hglue 4mm and 
 Elisabeth  Remm\thanks{Supported by the Laboratoire de Math. et Applications, Mulhouse and ASCR, Praha.}%
 } 

 \maketitle 
 \baselineskip 15pt plus 2pt minus 1 pt 

\begin{abstract} 
We investigate algebras with one operation.
We study when these algebras form a monoidal category
and analyze Koszulness and cyclicity of the corresponding operads. We also
introduce a new kind of symmetry for operads, the {\em dihedrality\/},
responsible for the existence of dihedral cohomology.

The main trick, which we call the {\em
polarization\/}, will be to represent an algebra with
one operation without any specific symmetry as an algebra with one
commutative and one anticommutative operations. We will 
try to convince the reader that this change of perspective might 
sometimes lead to new insights and results.

This point of view was used in~\cite{livernet-loday:unpubl} to introduce a
one-parametric family of operads whose specialization at $0$ is the
operad for Poisson algebras, while at a generic point it equals
the operad for associative algebras. We study this family
and explain how it can be used to interpret the deformation quantization
($*$-product) in a neat and elegant way.
\end{abstract}

 \parskip0pt 
 \vskip 3mm 
 \noindent 
 {\bf Table of content:} \ref{sec:ex}. 
                    Some examples to warm up 
                      -- page~\pageref{sec:ex} 
 \hfill\break\noindent 
 \hphantom{{\bf Table of content:\hskip .5mm}} \ref{sec:monoidal-structures}. 
                      Monoidal structures 
                      -- page~\pageref{sec:monoidal-structures} 
 \hfill\break\noindent 
 \hphantom{{\bf Table of content:\hskip .5mm}}  \ref{sec:kozulness-cyclicity-dihedrality}. 
                      Koszulness, cyclicity and dihedrality 
                      -- page~\pageref{sec:kozulness-cyclicity-dihedrality}

 \parskip3pt 

 \section*{Introduction} 
 \label{sec:intro} 
 If not stated otherwise, all algebraic objects in this paper will be defined over a fixed field $\mathbb{K}$ of characteristic 0. 
 We assume basic knowledge of operads as it can be gained for example 
from~\cite{markl-shnider-stasheff:book},
 though the first two sections can be read without this knowledge. 
 We are going to study classes of {\em algebras with one operation\/} 
 $\cdot:V \otimes V \to V$ and axioms given as linear combinations of 
 terms of the form $v_{\sigma(1)}\cdot (v_{\sigma(2)} \cdot 
 v_{\sigma(3)})$ and/or $(v_{\sigma(1)}\cdot v_{\sigma(2)}) \cdot 
 v_{\sigma(3)}$, where $\sigma \in \Sigma_3$ is a permutation. 
All classical examples of algebras, such 
 as associative, commutative, Lie and, quite surprisingly, Poisson 
 algebras, are of this type. 

 In a fancier, operadic, language this means that we are going to
 consider algebras over operads $\calP$ of the form $\calP =
 \Gamma(E)/(R)$, where $\Gamma(E)$ denotes the free operad generated
 by a $\Sigma_2$-module $E$ placed in arity $2$ and $(R)$ is the ideal
 generated by a $\Sigma_3$-invariant subspace $R$ of $\Gamma(E)(3)$.
 Operads of this form are called {\em
 quadratic\/}~\cite[Definition~3.31]{markl-shnider-stasheff:book}.  We
 moreover assume that the $\Sigma_2$-module $E$ is {\em generated by
 one element\/}. This gives a precise meaning to what we mean by of an
 ``algebra with one operation.''

 It follows from an elementary representation theory that, as a 
 $\Sigma_2$-module, 
 \[ 
 \mbox {either } 
 (1) \hskip 2mm E\cong \id_2,\ \mbox { or } 
 (2) \hskip 2mm E \cong \sgn_2,\ \mbox {or}\ 
 (3) \hskip 2mm E \cong \mathbb{K}[\Sigma_2], 
 \] 
 where $\id_2$ is the one-dimensional trivial representation, $\sgn_2$ is 
 the one-dimensional signum representation and $\mathbb{K}[\Sigma_2]$ is the 
 two-dimensional regular representation. In terms of the product these 
 cases can be characterized by saying that $\cdot: V \otimes V \to V$ is 

 (1) commutative, that is $x \cdot y = y \cdot x$ for all $x,y \in V$, 

 (2) anticommutative, that is $x \cdot y = - y \cdot x$ for all $x,y \in V$, 

 (3) without any symmetry, which means that there is no relation 
 between $x \cdot y$ and $y \cdot x$.   

 \vskip 2mm 

 Any multiplication $\cdot : V 
 \otimes V \to V$ of type~(3) can decomposed into the sum of a 
 commutative multiplication $\bullet$ and an anti-commutative one 
 $[-,-]$ via the {\em polarization \/} given by 
 \begin{equation}
\label{pol}
 x \bullet y :=\frac{1}{\sqrt{2}}(x \cdot y+y \cdot x)\ 
 \mbox { and }\ 
 \left[x, y\right] :=\frac{1}{\sqrt{2}}(x \cdot y-y \cdot x),\ 
 \mbox { for }\ x,y \in V. 
 \end{equation}
 The inverse process of {\em depolarization\/} assembles a type~(1) 
 multiplication $\bullet$ with a type~(2) multiplication $[-,-]$ into 
 \begin{equation}
\label{depol}
 x \cdot y :=\frac{1}{\sqrt{2}}(x \bullet y+[x,y]),\ 
 \mbox { for }\ x,y \in V. 
 \end{equation}
 The coefficient $\frac{1}{\sqrt2}$ was chosen so that the polarization followed by the depolarization (and vice versa)
is the identity.
 On the operadic level, the above procedure 
 reflects the decomposition 
 \[ 
 \mathbb{K}[\Sigma_2] \cong \id_2 \oplus \sgn_2 
 \] 
 of the regular representation into the trivial and signum representations. 

 The polarization enables one to view structures with a type~(3) 
 multiplication (such as associative algebras in Example~\ref{ex:1})  as 
 structures with one commutative and one anticommutative operation, 
 while the depolarization interprets structures with one commutative 
 and one anticommutative operation (such as Poisson algebras in 
 Example~\ref{ex:poiss}) as structures with one type~(3) operation. We will 
 try to convince the reader that this change of perspective might 
 sometimes lead to new insights and results.  The 
 polarization-depolarization trick was independently employed by 
 Livernet and Loday in their unpublished preprint 
\cite{livernet-loday:unpubl}. 

\noindent
{\bf Acknowledgment.} We would like to express our thanks
to Jean-Louis Loday who gave us the unpublished 
manuscript~\cite{livernet-loday:unpubl} and to Petr Somberg for
teaching us some basic facts from representation theory which we
should have remembered from our respective kindergardens.

 \section{Some examples to warm up} 
 \label{sec:ex} 

 In this section we give a couple of examples to illustrate the 
 (de)polarization trick.  We will usually omit the $\bullet$ 
 denoting a commutative multiplication and write simply $xy$ instead of 
 $x \bullet y$. 
   
 \begin{example} 
 \label{ex:1} 
 {\rm 
 {\em Associative algebras\/} are traditionally understood as 
 structures with one operation of type~(3). Let us see what happens with the 
 associativity 
 \begin{equation} 
 \label{assoc} 
 (x \cdot y) \cdot z=x \cdot (y \cdot z), \ \mbox { for } 
 x,y,z \in V, 
 \end{equation} 
 if we polarize the multiplication $\cdot : V \ot V \to V$. 
 We claim that decomposing $x \cdot y=\frac{1}{\sqrt{2}}(xy+[x,y])$, the associativity~(\ref{assoc}) becomes
 equivalent to the following two axioms: 
 \begin{eqnarray} 
 \label{eq:3} 
 [x,y z]&=& [x,y] z +y[x,z], 
 \\ 
 \label{eq:1} 
 [y,[x,z]] &=& (x y) z-x (y  z). 
 \end{eqnarray}
 To verify this, observe that (\ref{assoc})  implies
 $$
 \begin{array}{c}
 u_1(x,y,z):=A(x,y,z)-A(y,x,z)+A(z,y,x)+A(x,z,y)+A(y,z,x)-A(z,x,y) =  0,\\
 u_2(x,y,z):=A(x,y,z)+A(y,x,z)-A(z,y,x)+A(x,z,y)-A(y,z,x)-A(z,x,y) =0,
 \end{array}
 $$
 with 
\begin{equation} 
\label{A} 
A(x,y,z):=(x\cdot y) \cdot z-x\cdot (y\cdot z)
\end{equation}
the associator. Now, axiom~(\ref{eq:3}) (resp.~(\ref{eq:1})) 
can be obtained from
$u_1(x,y,z)=0$ (resp.~from $u_2(x,y,z)=0$) 
by the depolarization, that is replacing $A(x,y,z)$ by  
$$\frac{1}{2}\left\{(xy)z+[xy,z]+ [x,y]z+[[x,y],z]-x(yz)-x[y,z]-[x,yz]-[x,[y,z]]\right\}. $$
To prove that (\ref{eq:3}) together with (\ref{eq:1}) 
imply  (\ref{assoc}), observe that   
$$A(x,y,z)=\frac{1}{4}\left\{(u_1+u_2)(x,z,y)-(u_1-u_2)(z,x,y)\right\}.$$
Let us remark that the summation of~(\ref{eq:1}) over cyclic 
permutations gives the 
 Jacobi identity 
 \begin{equation} 
 \label{eq:2} 
 [x ,[ y , z]] + [y ,[ z , x]]+[z ,[ x , y]]=0. 
 \end{equation} 
 } 
 \end{example}

 \begin{example} 
 \label{ex:poiss} 
 {\rm 
{\em Poisson algebras\/} are usually defined as structures with two 
 operations, a commutative associative one and an anti-commutative 
 one satisfying the Jacobi identity. These operations are 
 tied up by a distributive law 
 \[ 
 [x,y z]=[x,y]z +y [x,z] 
 \] 
 which we already saw in~(\ref{eq:3}). The depolarization reinterprets 
 Poisson algebras as structures with one type~(3) operation $\cdot: V 
 \ot V \to V$ and one axiom: 
 \begin{equation}
 \label{poiss} 
 x \cdot (y \cdot z)=(x \cdot y)\cdot z -\frac{1}{3} 
 \left\{(x \cdot z)\cdot y+(y \cdot z)\cdot x-(y \cdot 
 x)\cdot z-(z \cdot x)\cdot y)\right\}. 
 \end{equation} 
 To see this, write each of the expressions 
$[x,y z],\ [x,y]z,  \ [x,[y ,z]]$ and $x(y z)$
 as a linear sum of permutations of the expressions 
 $x \cdot (y \cdot z)$ and $(x \cdot y) \cdot z$ as follows:
 \begin{eqnarray}
 \nonumber 
 [[x,y],z]=\frac{1}{2}\left\{(x \cdot y - y \cdot x) 
\cdot z-z \cdot(x \cdot y-y \cdot x)\right\},\\
 \nonumber 
 (xy)z=\frac{1}{2}\left\{(x \cdot y + y \cdot x) \cdot z+z 
\cdot(x \cdot y+y \cdot x)\right\},\\
 \nonumber
 \left[x,y\right]z=\frac{1}{2}\left\{(x \cdot y - y \cdot x) \cdot z+z \cdot
(x \cdot y-y \cdot x)\right\},\\
 \nonumber
 \left[xy,z\right]=\frac{1}{2}\left\{(x \cdot y + y \cdot x) \cdot z-z \cdot
(x \cdot y+y \cdot x)\right\}.
 \end{eqnarray}
Then associativity  of the commutative part becomes
 $$
\begin{array}{lll}
v_1(x,y,z) & := & (x \cdot y) \cdot z+(y \cdot x) 
\cdot z-(z \cdot y) \cdot x -(y \cdot z) \cdot x \\
& & -x \cdot (y \cdot z)-x \cdot (z \cdot y)+z 
\cdot (y \cdot x)+z \cdot (x \cdot y)=0.
\end{array}
$$
The Jacobi identity gives
$$
\begin{array}{ll} 
v_2(x,y,z) := & (x \cdot y) \cdot z-(y \cdot x) \cdot z-(z \cdot y) \cdot x-(x \cdot z) \cdot y+(y \cdot z) \cdot x+(z \cdot x) \cdot y \\
& -x\cdot (y \cdot z)+y \cdot (x \cdot z)+z \cdot (y \cdot x)+x \cdot (z \cdot y)-y \cdot (z \cdot x)-z \cdot (x \cdot y)=0.
\end{array}
$$
Finally, the distributive law gives
$$
\begin{array}{lll} 
v_3(x,y,z) & := & (x \cdot y)\cdot z-(y \cdot x) \cdot z+(z \cdot y) \cdot x+(x \cdot z) \cdot y \\ 
& & +(y \cdot z) \cdot x-(z \cdot x) \cdot y -x \cdot (y \cdot z)+y \cdot (x \cdot z)\\
& & - z \cdot (y \cdot x)-x \cdot (z \cdot y)-y \cdot (z \cdot x)+z \cdot (x \cdot y)=0.
\end{array}
$$
Now it is straightforward to verify that the vector
 $$v(x,y,z):=(x \cdot y)\cdot z-x \cdot (y \cdot z)
-\frac{1}{3}\left\{(x \cdot z) \cdot y+(y \cdot z) \cdot y-(y \cdot x)\cdot z-(z \cdot x) \cdot y\right\}.$$
determined by~(\ref{poiss}) can be expressed as
\begin{equation}
\label{eq}
v(x,y,z)=\frac{1}{6}\left\{2v_1(x,y,z)+v_2(x,y,z)+v_3(x,y,z)+2 v_3(z,x,y)\right\}.
\end{equation}
 This shows that the depolarized multiplication indeed fulfills~(\ref{poiss}). To prove that~(\ref{poiss}) implies
the usual axioms of Poisson algebras, one might use a similar 
straightforward method as in Example~\ref{ex:1}. Operadically this 
means the following. The operad $\Poiss$ for 
Poisson algebras can be presented as 
\[
\Poiss=\Gamma (\mathbb{K}[\Sigma_2])/(R),
\]
where $R$ is the 6-dimensional 
$\Sigma_3$-invariant subspace of  $\Gamma (\mathbb{K}[\Sigma_2])(3)$ 
generated by $v_1$
(the associativity), $v_2$ (the Jacobi) and $v_3$ (the distributive law). Equation~(\ref{eq}) implies that 
$v \in R$ but we must prove that $v$ in 
fact generates $R$. We leave this as an exercise to the reader.

The depolarization form of the axioms 
was already used to study deformations and rigidity of Poisson algebras in 
\cite{goze-remm:in_preparation}.
}
 \end{example}

 \begin{example} 
 {\rm 
 In this example, the ground field will be the complex numbers $\mathbb{C}$. In their unpublished 
 note  \cite{livernet-loday:unpubl}, Livernet and Loday  
  considered a one-parametric family of algebras with the axioms 
 \begin{eqnarray} 
 \nonumber 
 [x ,[ y , z]] + [y ,[ z , x]]+[z ,[ x , y]] &=& 0, 
 \\ 
 \label{eq:6} 
 [x,y z]&=& [x,y] z +y[x,z], 
 \\ 
 \label{eq:11}
  (x y) z-x (y  z) &=& q[y,[x,z]], 
 \end{eqnarray} 
 depending on a complex parameter $q$.  Observe that, for 
 $q \neq 0$, the first axiom (the Jacobi identity) is implied by the 
 third one. Let us call algebras satisfying the above axioms {\em 
 $LL_q$-algebras\/} (from Livernet-Loday). 

 For $q=0,$ ~(\ref{eq:11}) becomes the associativity and  we recognize the usual definition of Poisson algebras. 
 If $q=1$, we get associative algebras, in the polarized form of 
 Example~\ref{ex:1}. Furthermore, one may also consider the limit for $q \rightarrow \infty$: 
 \begin{eqnarray} 
 \nonumber 
 [x ,[ y , z]] + [y ,[ z , x]]+[z ,[ x , y]] &=& 0, 
 \\ 
 \nonumber 
 {}[x,y z]&=& [x,y] z +y[x,z], 
 \\ 
 \nonumber 
 [y,[x,z]] &=& 0. 
 \end{eqnarray} 
 In this case, the first identity trivially follows from the last 
 one.  These $LL_\infty$-algebras are algebras with a two-step-nilpotent 
 anticommutative bracket and a commutative multiplication, related by a 
 distributive law (the second equation). 

 The depolarization allows one to interpret $LL_q$-algebras as algebras 
 with one type~(3) operation $\cdot : V \ot V \to V$. The corresponding 
 calculation was made in~\cite{livernet-loday:unpubl}. One must 
 distinguish two cases. For $q \neq -3$ we get the axiom 
 \[ 
 x \cdot (y \cdot z)=(x \cdot y)\cdot z +\frac{q-1}{q+3}
\left\{ 
 (x \cdot z)\cdot y+(y \cdot z)\cdot x-(y \cdot 
 x)\cdot z-(z \cdot x)\cdot y)\right\}, 
 \] 
 while for $q = -3$ we get a structure with three axioms 
 \begin{eqnarray*} 
(x \cdot z)\cdot y+(y \cdot z)\cdot x-(y \cdot 
 x)\cdot z-(z \cdot x)\cdot y)&=&0, 
 \\ 
 A(x,y,z)+A(z,y,x)&=&0, 
 \\ 
 A(x,y,z)+A(y,z,x)+A(z,x,y)&=&0, 
 \end{eqnarray*} 
where $A$ denotes, as usual, the 
associator~(\ref{A}). It can be easily verified that the formula
\begin{equation}
\label{za_chvili_prednaska}
x \star y:= \frac{1+\sqrt{q}}{2}x \cdot y +\frac{1-\sqrt{q}}{2} y
\cdot x
\end{equation}
converts $LL_q$-algebras, for $q\neq 0, \infty$, into associative algebras.
Operadically this means that the operad $\mathcal{LL}_q$ 
for $LL_q$-algebras is, for $q \not\in \{0,\infty\}$, 
isomorphic to $\mathcal{LL}_1=\mathcal{A}ss$. This fact was observed also 
in ~\cite{livernet-loday:unpubl}. 
} 
 \end{example} 

\begin{example} 
 {\rm\ 
In this example we explain how Livernet-Loday 
algebras can be used to interpret deformation quantization 
of Poisson algebras. The ground ring here will be the ring 
${\KK}[[t]]$ of formal power series in $t$.
Let us recall~\cite{bayen78} that a {\em $*$-product\/} on a 
$\KK$-vector space $A$ is a $\KK[[t]]$-linear associative unital 
multiplication
$* : A[[t]] \ot A[[t]] \to A[[t]]$ 
which is commutative mod~$t$. Expanding, for $u,v \in A$, 
\begin{equation}
\label{kacer}
u*v = u*_0v + t\ u*_1v  + t^2\  u*_2v + \cdots,\ \mbox { with } u*_iv \in A
\mbox { for } i \geq 0, 
\end{equation}
one easily verifies that the operations $\cdot_0$ and $[-,-]_0$
defined by
\begin{equation}
\label{Raleigh}
u \cdot_0 v : =  u*_0v
\ \mbox{ and }\ [u,v]_0 :=  u*_1v -  v*_1u, \ u,v \in A,
\end{equation}
are such that $P:= (A,\cdot_0,[-,-]_0)$ is a Poisson algebra. 
The object $(A[[t]],*)$ is sometimes also called the {\em 
deformation quantization\/} of the Poisson algebra $P$. In applications, $P$ is the $\mathbb{R}$-algebra  
$C^\infty(M)$ of smooth functions on a Poisson manifold $M$ that represents the phase space of a classical physical system.
One moreover assumes that all products $*_i$, $i \geq 0$, 
in~(\ref{kacer}) are bilinear differential operators, see again~\cite{bayen78} 
for details. The relevance of $LL_q$-algebras for
quantization is explained in the following theorem.
 
\begin{theorem}
A $*$-product on a $\mathbb{K}$-vector space $A$ is 
the same as an $LL_{t^2}$-algebra structure on the $\KK[[t]]$-module $V : =A[[t]]$.  
\end{theorem}

\noindent{\bf Proof.}
Given a $*$-product,  
define $\bullet : V \ot V \to V$ and $[-,-] : V \ot V \to V$ as the 
polarization~(\ref{pol}) of $* : V \ot V \to V$. Commutativity of $*$ mod $t$ 
means that $[-,-] = 0$ mod $t$, therefore there exists a bilinear antisymmetric
map $\{-,-\} : V \ot V \to V$ such that $[-,-] = t\{-,-\}$. It is immediate 
to check that $(V,\bullet ,\{-,-\})$ is an $LL_{t^2}$-algebra.
On the other hand, given an $LL_{t^2}$-algebra $(V,\bullet ,\{-,-\})$, then 
\[
u*v := \frac{1}{\sqrt{2}}(u \bullet v+t \{u,v\}),\ 
 \mbox { for }\ u,v \in V, 
\]
clearly defines a $*$-product on $A$.%
\qed
}
\end{example}

 \begin{example}
\label{Lie-adm} 
 {\rm\ 
 Recall~\cite{goze-remm:ja04} that a type~(3) product $\cdot:V \otimes V 
 \rightarrow V$ is called {\em Lie-admissible\/} if the commutator 
 of this product is a Lie bracket or, equivalently, 
 that the antisymmetric part $[-,-]$ of its polarization fulfills the 
 Jacobi identity~(\ref{eq:2}). This observation suggests that the polarization might be particularly 
 suited for various types of Lie-admissible algebras. 

 Some  important classes of Lie-admissible algebras where studied 
 in~\cite{goze-remm:2}. Before we recall the definitions,
  we note that a type~(3) product $\cdot : V \ot V \to V$ is 
 Lie-admissible if and only if its associator~(\ref{A}) satisfies 
 \begin{equation} 
 \label{eq:5} 
 \sum_{\sigma \in \Sigma_3}(-1)^{\epsilon(\sigma)} 
 A(x_{\sigma(1)},x_{\sigma(2)},x_{\sigma(3)}) = 0, 
 \end{equation} 
 where $\epsilon(\sigma)$ denotes the signature of the permutation 
 $\sigma$. Now {\em $G$-associative algebras\/},  
 where $G$ is a (not necessary 
 normal) subgroup of $\Sigma_3$, are algebras with a type~(3) 
 multiplication whose commutator satisfies a condition which is, for $G \neq \Sigma_3,$
 stronger  
 than~(\ref{eq:5}), namely 
 \[ 
 \sum_{\sigma \in G}(-1)^{\epsilon(\sigma)} 
 A(x_{\sigma(1)},x_{\sigma(2)},x_{\sigma(3)}) = 0. 
 \] 
Therefore we have six different types of $G$-associative algebras 
 corresponding to the following six subgroups of~$\Sigma_3$: 
 \[ 
 G_1 : =\{1\},\  G_2 :=\{1,\tau_{12}\}, 
 \ G_3 := 
\{1,\tau_{23}\},\ G_4:=\{1,\tau_{13}\},
\ G_5:=A_3  \mbox { and } G_6 :=\Sigma_3, 
\] 
where $\tau_{ij}$ denotes the transposition $i \leftrightarrow j$ and $A_3$ is 
 the alternating subgroup of~$\Sigma_3$. 
 {\em $G_1$-associative algebras\/} are clearly associative algebras whose polarization we 
 discussed in Example~\ref{ex:1}. 
 {\em $G_2$-associative algebras\/} are the same as 
 {\em Vinberg algebras\/}, also called {\em left-symmetric
 algebras\/},~see~\cite{vinberg:trudy63}. 

 In the polarized form, 
$G_2$-associative algebras are structures with a commutative 
 multiplication  and a Lie bracket related by the axiom: 
 \begin{equation} 
\label{equ:G2} 
2[x, y] z+[[x ,y], z]-x (y z)+y(x z) -x [y,z]+y [x,z]-[x,y z]+[y,x
 z]=0.  \end{equation} As suggested by Loday,~(\ref{equ:G2}) can be
 written as the sum
\[
\left\{\rule{0em}{1em} \hskip -.1em
(xz)y -x(zy) - [z,[x,y]] 
\hskip -.1em\right\}
+
\left\{\rule{0em}{1em}\hskip -.1em
[x,y]z + y[x,z] -[x,yz]
\hskip -.1em\right\} +
\left\{\rule{0em}{1em}\hskip -.1em
[y,xz] - x[y,z] - [y,x]z
\hskip -.1em\right\} = 0
\]
of three terms which vanish separately if the multiplication is
associative, see~(\ref{eq:3}) and~(\ref{eq:1}).
 {\em $G_3$-associative algebras\/} are also classical objects, 
 known as {\em right-symmetric algebras}~\cite{matsushima:OJMa68} or 
{\em pre-Lie algebras\/}~\cite{gerstenhaber:AM63}.
 
The operad $\mathcal{P}re$-$\mathcal{L}ie$ associated to pre-Lie 
algebras is isomorphic to the operad $\Vinb$ 
for Vinberg algebras via the operadic homomorphism determined by
$$
\begin{array}{l}
xy \mapsto  xy, \quad \left[x,y\right] \mapsto - \left[x,y\right].
\end{array}
$$
This homomorphism converts~(\ref{equ:G2}) into
\begin{eqnarray*} 
  2[x, y] z-[[x ,y], z]+x (y z)-y(x z) -x [y,z]+y [x,z]-[x,y z]+[y,x z]=0.& 
 \end{eqnarray*} 

 After polarizing, we identify {\em $G_4$-associative algebras\/} as structures 
 satisfying 
 $$ 
 (x y) z-x (y z)=[[x,z],y], 
 $$ 
 which clearly implies the Jacobi~(\ref{eq:2}). 
 {\em $G_5$-associative algebras} 
 have a commutative multiplication and a Lie bracket 
 tied together by
 \begin{equation} 
 \label{eq:7} 
 [x  y,z] +[y  z, x ]+[z x,y]=0. 
 \end{equation}

 The polarization of {\em $G_6$-associative algebras\/}, which are sometimes confusingly called 
 just {\em Lie-admissible\/} algebras, reveals that the category of 
 these objects is a dull one, consisting of structures with a 
 commutative multiplication and a Lie bracket, with no relation between 
 these two operations. 

 Let us close this example by observing
that axiom~(\ref{eq:6}) of $LL_q$-algebras implies 
 axiom~(\ref{eq:7}) of $G_5$-associative algebras, therefore $LL_q$-algebras form, 
 for each $q$, a subcategory of the category of $G_5$-associative algebras. 
 An equally simple observation is that the polarized product of $G_4$-associative algebras satisfies 
the first and the third identities of $LL_{-1}$-algebras but not
the distributive law. 
} 
 \end{example} 

\begin{example}
{\rm\ Lie-admissible structures mentioned in Example~\ref{Lie-adm} are
rather important. As it was shown in the seminal
paper~\cite{gerstenhaber:AM63}, there exists a natural pre-Lie
structure on the Hochschild cochain complex of every associative
algebra induced by an even more elementary structure christened later,
in~\cite{gerstenhaber-voronov:IMRN95}, a {\em brace algebra\/}. This
pre-Lie structure is responsible for the existence of the {\em
intrinsic bracket\/} on the Hochschild cohomology, see
again~\cite{gerstenhaber:AM63}.

We offer the following generalization of this structure.
For a vector space $V$, denote by
\[
X := 
\bigoplus_{m,n \geq 1}\Lin (\otexp Vm,\otexp Vn)
\]
the space of all multilinear maps.
{}For $f \in \Lin (\otexp Vb,\otexp Va)$ and $g \in \Lin (\otexp
Vd,\otexp Vc)$ define $f \circ^i_j g \in \Lin (\otexp V{a +
c-1},\otexp V{b + d -1})$ to be the map obtained by composing the
$j$-th output of $g$ into the $i$-th input of $f$ and arranging the
remaining outputs and inputs as indicated in Figure~\ref{gee!}.
\begin{figure}
\begin{center}
\unitlength 1.2mm 
\begin{picture}(32,30)(-30,5)
\thicklines
\put(-13,13){\line(0,-1){5}}
\put(-13,23){\line(0,1){5}}
\put(-13,28){\line(-1,0){11}}
\put(-24,23){\line(1,0){11}}
\put(-18,23){\line(0,-1){10}}
\bezier{211}(-22,23)(-22,18.5)(-26,16)
\bezier{211}(-26,16)(-30,13.5)(-30,9)
\bezier{211}(-20,23)(-20,18.5)(-24,16)
\bezier{211}(-24,16)(-28,13.5)(-28,9)
\bezier{211}(-24,23)(-24,18.5)(-28,16)
\bezier{211}(-28,16)(-32,13.5)(-32,9)
\put(-24,23){\line(-1,0){1}}
\put(-25,23){\line(0,1){5}}
\put(-25,28){\line(1,0){1}}
\bezier{211}(-16,23)(-16,18.5)(-12,16)
\bezier{211}(-14,23)(-14,18.5)(-10,16)
\bezier{211}(-10,16)(-6,13.5)(-6,9)
\bezier{211}(-12,16)(-8,13.5)(-8,9)
\bezier{211}(-20,13)(-20,17.5)(-24,20)
\bezier{211}(-22,13)(-22,17.5)(-26,20)
\bezier{211}(-24,20)(-28,22.5)(-28,27)
\bezier{211}(-26,20)(-30,22.5)(-30,27)
\bezier{211}(-16,13)(-16,17.5)(-12,20)
\bezier{211}(-14,13)(-14,17.5)(-10,20)
\bezier{211}(-12,20)(-8,22.5)(-8,27)
\bezier{211}(-10,20)(-6,22.5)(-6,27)
\put(-13,13){\line(-1,0){10}}
\put(-23,13){\line(0,-1){5}}
\put(-23,8){\line(1,0){10}}
\put(-28,27){\line(0,1){5}}
\put(-30,27){\line(0,1){5}}
\put(-8,27){\line(0,1){5}}
\put(-6,27){\line(0,1){5}}
\put(-15,28){\line(0,1){4}}
\put(-19,28){\line(0,1){4}}
\put(-23,28){\line(0,1){4}}
\put(-6,9){\line(0,-1){5}}
\put(-8,9){\line(0,-1){5}}
\put(-28,9){\line(0,-1){5}}
\put(-30,9){\line(0,-1){5}}
\put(-32,9){\line(0,-1){5}}
\put(-14,8){\line(0,-1){4}}
\put(-16,8){\line(0,-1){4}}
\put(-18,8){\line(0,-1){4}}
\put(-20,8){\line(0,-1){4}}
\put(-22,8){\line(0,-1){4}}
\put(-6,4){\vector(0,1){2}}
\put(-8,4){\vector(0,1){2}}
\put(-14,4){\vector(0,1){2}}
\put(-16,4){\vector(0,1){2}}
\put(-18,4){\vector(0,1){2}}
\put(-20,4){\vector(0,1){2}}
\put(-22,4){\vector(0,1){2}}
\put(-28,4){\vector(0,1){2}}
\put(-30,4){\vector(0,1){2}}
\put(-32,4){\vector(0,1){2}}
\put(-6,30){\vector(0,1){2}}
\put(-8,30){\vector(0,1){2}}
\put(-15,30){\vector(0,1){2}}
\put(-19,30){\vector(0,1){2}}
\put(-23,30){\vector(0,1){2}}
\put(-28,30){\vector(0,1){2}}
\put(-30,30){\vector(0,1){2}}
\put(-19,25.5){\makebox(0,0)[cc]{$f$}}
\put(-18,10){\makebox(0,0)[cc]{$g$}}
\end{picture}
\end{center}
\label{gee!}
\caption{The composition $f \circ^4_3 g \in \Lin (\otexp V{10},\otexp
V7)$ of functions $f \in \Lin (\otexp V6,\otexp V3)$ and $g \in \Lin
(\otexp V5,\otexp V5)$.}
\end{figure}
Define finally
\[
f \circ g := \sum_{{1 \leq i \leq b},{\ 1 \leq j \leq c}}
             (-1)^{i(b+1) + j(c+1)}
f \circ^i_j g.
\]
We leave to the reader to verify that $(X,\circ)$ {\em is a $G_6$-associative
algebra\/}. 

Let $\Yhoch \subset X$ be the subspace $\Yhoch : = \bigoplus_{m \geq
1}\Lin (\otexp Vm,V)$ and let dually $\Ycohoch := \bigoplus_{n \geq
1}\Lin (V,\otexp Vn)$. Clearly both $\Yhoch$ and $\Ycohoch$ are
$\circ$-closed.  It turns out that $(\Yhoch,\circ)$ is a {\em
$G_3$-associative algebra\/} (= pre-Lie algebra) and
$(\Ycohoch,\circ)$ a {\em $G_2$-associative algebra\/} (= Vinberg
algebra).  We recognize $(\Yhoch,\circ)$ as the underlying space of
the Hochschild cochain complex $C^*_{\it Hoch}(A;A)$ of an associative
algebra $A = (V,\cdot)$ with the classical pre-Lie
structure~\cite{gerstenhaber:AM63}. The space $(\Ycohoch,\circ)$ has a
similar interpretation in terms of the Cartier cohomology of
coassociative coalgebras~\cite{cartier}.

To interpret $X$ in a similar way, wee need to recall that an {\em
infinitesimal bialgebra\/}~\cite{aguiar:02} (also called a {\em mock
bialgebra\/} in~\cite{fox-markl:ContM97}) is a triple $(V,\mu,\delta)$,
where $\mu$ is an associative multiplication, $\delta$ is a
coassociative comultiplication and
\[
\delta (\mu(u,v)) = \delta_{(1)}(u) \otimes \mu (\delta_{(2)}(u),v)
+ \mu(u, \delta_{(1)}(v)) \otimes \delta_{(2)}(v)
\]
for each $u,v \in V$, with the standard Sweedler's notation for the
comultiplication used.  It turns out that $X$ is the underlying space
of the cochain complex defining the cohomology of an infinitesimal
bialgebra and $[f,g] := f \circ g - g \circ f$ is the intrinsic
bracket, see~\cite{markl:ib}, on this cochain complex.  The reader is
encouraged to verify that the axioms of infinitesimal bialgebras can
be written as the `master equation'
\[
[\mu + \delta,\mu + \delta] = 0,
\] 
with $\mu : V \ot V \to V$ and $\delta : V \to V \ot V$ interpreted 
as elements of~$X$.
}
\end{example}

 \section{Monoidal structures} 
 \label{sec:monoidal-structures} 

 Consider the category $\Palg$ of algebras over a fixed operad 
 $\calP$. Following~\cite{getzler-jones:preprint} we say that $\calP$ 
 is a {\em Hopf operad\/}, if the category $\Palg$ admits a strict 
 monoidal structure $\odot : \Palg \times \Palg \to \Palg$ such that 
 the forgetful functor $\BOX : \Palg \to \Vect_\KK$ to the category 
 of $\KK$-vector spaces with the standard tensor product, 
 is a strict monoidal morphism, see~\cite[VII.1]{maclane:71} for the 
 terminology.  This condition can be expressed solely in terms of 
 $\calP$ as in the following definition.

 \begin{definition} 
 An operad $\mathcal{P}$ is a Hopf operad if there exists an operadic 
 map $\Delta:\mathcal{P} \rightarrow \mathcal{P} \otimes \mathcal{P}$ 
 (the {\em diagonal\/}) which is coassociative in the sense that 
 \begin{equation} 
 \label{coass} 
 (\Delta \otimes \id_\calP)\Delta =(\id_\calP \otimes \Delta)\Delta, 
 \end{equation} 
 where $\id_\calP : \calP \to \calP$ denotes the identity.  We also 
 assume the existence of a {\em counit\/} $e:\mathcal{P} \rightarrow 
 \Com$, where $\Com$ is the operad for commutative associative 
 algebras, satisfying 
 \begin{equation} 
 \label{dymka} 
 (e \otimes \id_\calP ) \Delta = (\id_\calP \otimes e) \Delta = \id_\calP. 
 \end{equation} 
 \end{definition} 
 
The last equation uses the canonical identification
 $\calP \cong \Com \ot \calP \cong \calP \ot \Com$.  Our terminology
 slightly differs from the original one
 of~\cite{getzler-jones:preprint} which did not assume the counit.
 Our assumption about the existence of the
 counit rules out trivial diagonals.
 
 The diagonal $\Delta : \mathcal{P} \rightarrow \mathcal{P} \otimes
 \mathcal{P}$ induces a product $\odot : \Palg \times \Palg \to \Palg$
 in a way described for example
 in~\cite[page~197]{markl-shnider-stasheff:book}. Equation~(\ref{coass})
 is equivalent to the coassociativity of this product. To
 interpret~(\ref{dymka}), observe that, since $\Com$ is
 isomorphic to the endomorphism operad $\End_\KK$ of the ground
 field, the counit $e$ equips $\KK$ with a $\calP$-algebra structure.
 Equation~(\ref{dymka}) then says that $\KK$ with this structure is
 the unit object for the monoidal structure induced by $\Delta$.

 In the rest of this section we want to discuss the existence Hopf 
 structures on quadratic operads $\calP = \Gamma(E)/(R)$ with one 
 operation. Let us look more closely at the map $e : \calP \to 
 \Com$ first. Since $\Com = \Gamma(\id_2)/(R_{\rm ass})$, with $\id_2$ the 
 trivial representation of $\Sigma_2$ and $(R_{\rm ass})$ the ideal 
 generated by the associativity, the counit $e$ is determined by a 
 $\Sigma_2$-equivariant map 
 \begin{equation} 
 \label{kourim_vodni_dymku} 
 e(2) : E \to \id_2. 
 \end{equation} 
 If $E = \id_2$ (case~(1) of the nomenclature of the introduction), such a map is the multiplication by a scalar 
 $\alpha$. If $E = \sgn_2$ (case~(2)), the only equivariant $e(2)$ is the 
 zero map. Finally, if $E  = \KK[\Sigma_2]$ (case~(3)), $e(2)$ must be the 
 projection $\KK[\Sigma_2] \to \id_2$ multiplied by some $\alpha \in \KK$. 
 
 Equation~(\ref{dymka}) implies the non-triviality of $e(2)$. This
 excludes case~(2) and implies that $\alpha \not= 0$ in cases~(1)
 and~(3). In these two cases we may moreover assume the {\em normalization\/}
 $\alpha = 1$, the general case can be brought to this form by rescaling
 $e \mapsto \alpha^{-1} e$, $\Delta \mapsto \alpha \Delta$.
 
 Let us introduce the following useful pictorial language.  Denote by
 $\tri {}{} \in \calP(2)$ the operadic generator for a type~(3)
 operation (a multiplication with no symmetry). Similarly, we denote
 the generator for a commutative operation by $\jedna {}{}$ and for an
 anti-commutative one by $\dva {}{}$. The right action of the
 generator $\tau \in \Sigma_2$ on $\calP(2)$ is, in this language, described by
\[
\tri 12 \tau  = \tri 21,\ \jedna 12 \tau = \jedna 21 = \jedna 12 
\mbox {\hskip 3mm\ and \hskip 3mm}
\dva 12 \tau = \dva 21 = -\dva 21.
\]
The polarization~(\ref{pol}) is then given by
 \[ 
 \jedna 12 =\frac{1}{\sqrt{2}}\left(\tri 12+ \tri 21\right)\ 
 \mbox { and \ }\ 
 \dva 12  :=\frac{1}{\sqrt{2}}\left(\tri 12- \tri 21\right),
 \] 
and the depolarization~(\ref{depol}) by
\[ 
\tri 12  :=\frac{1}{\sqrt{2}}\left(\jedna 12 + \dva 12   \right).
\]

In the rest of this section we investigate the existence of diagonals
for quadratic operads with one operation. Since the diagonal is, by
assumption, an operadic homomorphism, it is uniquely determined by its value on
a chosen generator of $\calP(2)$. Let us see what can be concluded from this
simple observation. As before, we distinguish three cases.

\vskip 3mm
\noindent 
{\em Case~(1).\/} In this case, the operad $\calP$ is generated by one
commutative bilinear operation $\jedna {}{} \in \calP(2)$. The
diagonal must necessarily satisfy
\[
\Delta(\jedna 12 )=A \left( \jedna 12\otimes \jedna 12\right), \mbox {
for some $A \in \KK$.}
\]
The coassociativity~(\ref{coass}) is  fulfilled automatically
while~(\ref{dymka}) implies $A=1$.

\vskip 3mm
\noindent 
{\em Case~(2).\/} Analyzing the counit, we already observed that
operads with one anti-symmetric operation do not admit a (counital)
diagonal. An easy argument shows that non-trivial diagonals for
type~(2) operads do not exists even if we do not demand the existence
of a counit. Indeed, in case~(2) we have $\calP(2) \cong \sgn_2$ while
$\calP(2) \otimes \calP(2) \cong \sgn_2 \ot \sgn_2 \cong \id_2$, therefore
$\Delta(2) : \calP(2) \to \calP(2) \ot \calP(2)$ is trivial, as is any
$\Sigma_2$ equivariant map $\sgn_2 \to \id_2$. Let
us formulate this observation as:

 \begin{theorem} 
\label{thm:1}
   There is no non-trivial diagonal on a quadratic operad generated by
   an anti-symmetric product. In particular, the operad $\Lie$ for Lie
   algebras is not an Hopf operad.
 \end{theorem}

\vskip 3mm
\noindent 
{\em Case~(3).\/} As an operadic homomorphism, the diagonal (if exists) is
uniquely determined by an element $\Delta(\tri 12) \in \calP(2) \otimes
\calP(2)$. The following proposition characterizes which choices of
$\Delta(\tri 12)$ may lead to a coassociative counital diagonal.

 \begin{proposition}
\label{prop:1} 
 Let $\mathcal{P}$ be a quadratic Hopf operad generated by a type~(3)
 product $\tri {}{}$.
Then there exists $B \in \KK$ such that the polarized form of the diagonal $\Delta$ is given by
 \begin{equation}
\label{eq:8}
 \Delta(\tri 12)= \tri 12 \otimes \tri 12  
  -B\left\{(\tri 12- \tri 21)\otimes(\tri 12- \tri 21) \right\}.
 \end{equation}
 The polarized version of this equation reads
\begin{equation}
\label{eq:9}
\left\{
  \begin{array}{rcl}
 \Delta(\jedna 12)&=& \displaystyle\frac{1}{\sqrt{2}} \left\{\jedna 12
 \otimes\jedna 12+ (1-4B)\left(\dva 12\otimes\dva 12 \right) \right\},
\\ 
 \rule {0pt}{23pt}\Delta(\dva 12)&=& 
\displaystyle\frac{1}{\sqrt{2}} \left\{\dva 12 
 \otimes\jedna 12 + \jedna 12\otimes\dva 12 \right\}. 
\end{array}
\right.
\end{equation}
 \end{proposition}

\vskip 2mm
\noindent
{\em Proof.} A simple bookkeeping.
  The most general choice for $\Delta(2): \calP(2) \to \calP(2) \ot
  \calP(2)$ is 
 \begin{equation}
\label{eq:10}
  \Delta(\tri 12)= A(\tri 12\otimes \tri 12) 
 +B(\tri 12 \otimes \tri 21)
  +C(\tri 21 \otimes \tri 12) 
 +D(\tri 21 \otimes \tri 21) 
 \end{equation}
 with \rule{0pt}{18pt}some $A,B,C,D \in \KK$.  
 A straightforward calculation shows that the coassociativity 
 $(\Delta \otimes \id_\calP)\Delta =(\id_\calP \otimes \Delta)\Delta$   
for $\Delta$ defined by~(\ref{eq:10}) has the following four families
 of solutions:
\[
\begin{array}{ll}
\mbox {(i)}\hskip 2mm D=C=0,\ A=B, &
\mbox {(ii)}\hskip 2mm B=C=-D,\ A \mbox { arbitrary},
\\
\mbox {(iii)}\hskip 2mm D=B=0,\ A=C, &\rule{0pt}{15pt}
\mbox {(iv)}\hskip 2mm B=C=D=A.
\end{array} 
\]
The counit condition~(\ref{dymka}) leads to the system:
\[
A+C=1,\ B+D=0,\  A+B=1,\ \mbox { and } C+D=0. 
\]
We easily conclude that the only solution is a type (ii) one with 
$B=C=-D$ and $A=1-B$. This gives~(\ref{eq:8}) whose polarization 
is~(\ref{eq:9}).%
\qed

Proposition~\ref{prop:1} offers a mighty tool to investigate the
existence of Hopf structures for type~(3) operads. It says that such an
operad $\calP$ is a Hopf operad if and only if there exists $B
\in \KK$ such that the diagonal defined by~(\ref{eq:8})
(resp.~(\ref{eq:9}) in the polarized form) extends to an operad map,
i.e. preserves the relations $R$ in the quadratic presentation
$\Gamma(E)/(R)$ of $\calP$.

Regarding the existence of diagonals in general, an operad
might admit no Hopf structure at all (examples of this situation are
provided by Theorem~\ref{thm:1}), it might admit exactly one Hopf
structure (see Example~\ref{ethm:2} for operads with this
property), or it might admit several different monoidal structures, as
illustrated in Example~\ref{ex:two}.

Let us formulate another simple proposition whose proof we leave as an
exercise.  We say that $\calP$ is a {\em set-operad\/}, if there
exists an operad ${\mathcal S}$ in the monoidal category of sets such
that, for any $n \geq 1$, $\calP(n)$ is the $\KK$-linear span of
${\mathcal S}(n)$, and  that the operad structure of
$\calP$ is naturally induced from the operad structure of ${\mathcal
  S}$.

\begin{proposition}
  Every set-operad $\calP$ admits an Hopf structure given by the
  formula $\Delta(p) := p \ot p$, for any $p \in \calP$.
\end{proposition}

 \begin{example} 
\label{ethm:2}
{\rm
Let $\mathcal{LL}_{q}$ denote the operad for $LL_q$ algebras. Then 
the one-parametric family $\left\{\mathcal{LL}_{q} \right\}_{q \neq \infty}$ is a family of Hopf operads, with the diagonal given by   
\begin{equation}
\label{eq:14}
 \Delta(\tri 12)= \tri 12 \otimes \tri 12  
  -\frac{1-q}{4}\left\{(\tri 12- \tri 21)\otimes(\tri 12- \tri 21) \right\}.
 \end{equation}
 The polarized version of this equation reads
\begin{equation}
\label{eq:15}
\left\{
  \begin{array}{rcl}
 \Delta(\jedna 12)&=&\displaystyle\frac{1}{\sqrt{2}} \left\{ (\jedna 12
 \otimes\jedna 12+ q\left(\dva 12\otimes\dva 12 \right) \right\},
\\ 
 \rule {0pt}{23pt}\Delta(\dva 12)&=& \displaystyle\frac{1}{\sqrt{2}} 
\left\{\dva 12 
 \otimes\jedna 12 + \jedna 12\otimes\dva 12\right\}. 
\end{array}
\right.
\end{equation}
The above normalized 
diagonal is moreover unique for each $q \neq \infty$. 
Observe that the limit for $q \to \infty$ of 
formulas~(\ref{eq:14}) (resp.~(\ref{eq:15})) 
does not make sense and, indeed, it can be easily shown that 
the operad $\mathcal{LL}_{\infty}$ is not Hopf.  
}\end{example}

\begin{example}
\label{ex:two}
{\rm We give a funny example of a category of algebras
  which admits several non-equivalent monoidal structures.  Let us
  consider a type~(3) product $x,y \mapsto x\cdot y$, with the axiom
 \[
 (x \cdot y) \cdot z 
 =z \cdot (y \cdot x).
\]
 Then 
  \[
 \Delta(\tri 12):= \tri 12 \otimes \tri 12  
  -B\left [(\tri 12- \tri 21)\otimes(\tri 12- \tri 21) \right].
 \]
 defines an Hopf structure for any $B \in \KK$. 

Above we saw an algebra
admitting a one-parametric family of non-equivalent monoidal structures. 
It would be interesting to see a structure that 
admits a {\em discrete\/} family of 
non-equivalent Hopf structures.
} 
\end{example}

\begin{example}
{\rm
It can be shown that the only $G$-admissible algebras that admit a
monoidal structure are associative algebras. In particular, neither
Vinberg nor pre-Lie algebras form a monoidal category.
}
\end{example}

 \section{Koszulness, cyclicity and dihedrality} 
 \label{sec:kozulness-cyclicity-dihedrality} 

\def\rada#1#2{{#1,\ldots,#2}}
\def\Sigmap{\Sigma^+}

In this section we study 
cyclicity~\cite{getzler-kapranov:CPLNGT95} of operads mentioned
in the previous sections.
We then introduce the notion of {\em dihedrality\/} of operads and
investigate this property. To complete the picture, we also list
results concerning Koszulness~\cite{ginzburg-kapranov:DMJ94} 
of some operads with one operation.

Let us recall first what is a cyclic operad.
Let $\Sigmap_n$ be the group of automorphisms of the set $\{\rada
0n\}$. This group is,
of course, isomorphic to the symmetric group $\Sigma_{n+1}$, but the
isomorphism is canonical only up to an identification
$\{\rada 0n\} \cong \{\rada 1{n+1}\}$. 
We interpret $\Sigma_n$
as the subgroup of $\Sigmap_n$ consisting of permutations
$\sigma \in \Sigmap_n$ with $\sigma(0)= 0$. If $\gamma^+_n
\in \Sigmap_n$ denotes the cycle $(\rada 0n)$, that is, the
permutation with $\gamma^+_n(0)=1,\
\gamma^+_n(1)=2,\ \ldots,\tau_n(n)=0$, then $\gamma^+_n$
and $\Sigma_n$ generate $\Sigmap_n$.

By definition, each operad
$\mathcal{P}$ has a natural right action of $\Sigma_n$ on 
each piece $\mathcal{P}(n)$, $n \geq 1.$ 
The operad $\mathcal{P}$ is {\em cyclic\/} if this action
extends, for any $n \geq 1$, to a $\Sigmap_n$-action in a 
way compatible with structure
operations. See~\cite[Definition~II.5.2]{markl-shnider-stasheff:book} 
or the original 
paper~\cite{getzler-kapranov:CPLNGT95} for a precise definition.

We already recalled in the introduction that an (ordinary) 
operad $\calP$ is quadratic if it can be presented as $\calP = \Gamma(E)/(R)$, 
where $E=\calP(2)$ and $R \subset \Gamma(E)(3)$. 
The action of $\Sigma_2$ on
$E$ extends to an action of $\Sigma^+_2$, via the sign
representation
$\sgn : \Sigma^+_2 \to \{ \pm 1\} \cong \Sigma_2$. It can be easily verified
that this action induces a cyclic operad structure on the free operad
$\Gamma(E)$. In particular, $\Gamma(E)(3)$ is a right
$\Sigma^+_3$-module. 
An operad $\calP$ as above is called {\em cyclic quadratic\/} 
if the space of
relations $R$ is invariant under the action of $\Sigma^+_3$. Since $R$
is, by definition, $\Sigma_3$-invariant, $\calP$ is cyclic quadratic 
if and only if $R$ is preserved by the action of the generator~$\gamma^+_3$.

\begin{remark}
\label{pisi_na_letisti}
{\rm
There are operads that are both quadratic and cyclic 
but not cyclic quadratic. The simplest example of this
exotic phenomenon is provided by the free operad $\Gamma(V_{2,2})$ generated
by the $2$-dimensional irreducible representation $V_{2,2}$ of $\Sigma_3 \cong
\Sigmap_2$ placed in arity $2$.
In general, an operad $\calP$ is cyclic quadratic if and only if it is both
quadratic and cyclic and if the $\Sigmap_2$-action on $\calP(2)$ is
induced from the operadic $\Sigma_2$-action on $\calP(2)$ via the
homomorphism $\sgn : \Sigma^+_2 \to \{ \pm 1\} \cong \Sigma_2$.
}
\end{remark}

Let us turn our attention to the cyclicity of operads for algebras
with one operation. Since, as proved 
in~\cite[Proposition~3.6]{getzler-kapranov:CPLNGT95}, each
quadratic operad with one operation of type~(1) or~(2) is cyclic
quadratic, we shall
focus on operads with a type (3) multiplication. 
The right action of the generator $\gamma^+_3 \in \Sigmap_3$ on
$\Gamma(\KK[\Sigma_2])(3)$ is described in the following table: 
\[ 
\def\cd{\cdot}\def\arraystretch{1.2}
 \begin{array}{rclrcl} 
((x \cdot y) \cdot z)\gamma^+_3 &=& x \cd (y \cd z), 
&  (x \cd (y \cd z))\gamma^+_3 &=& (x \cd y) \cd z, 
\\
 ((y \cd z) \cd x)\gamma^+_3 &=& (y \cd x) \cd z,
&  (y \cd (z \cd x))\gamma^+_3   &=& y \cd (x \cd z),
\\
((z\cd x) \cd  y)\gamma^+_3 &=& y \cd(z \cd x),
&(z \cd (x \cd y))\gamma^+_3  &=& (y \cd z) \cd x,
\\
((y \cd x) \cd z )\gamma^+_3 &=& x \cd (z \cd y),
&  (y \cd (x \cd z))\gamma^+_3 &=& (x \cd z) \cd y,
\\
((z \cd y) \cd x)\gamma^+_3 &=& z \cd (y \cd x),
&(z \cd(y \cd x ))\gamma^+_3 &=& (z \cd y) \cd x,
\\
((x \cd z)\cd y)\gamma^+_3 &=& (z \cd x) \cd y,
&
(x\cd (z \cd y))\gamma^+_3 &=& z \cd (x \cd y).
 \end{array} 
\] 
Using this table, it is easy to investigate the cyclicity of operads
with one operation, see Example~\ref{nespravedliva_odmena} where the
corresponding analysis was done for $G_5$-associative algebras. The
results are summarized in Figure~\ref{tab}.

\begin{figure}[t]
\[
\def\arraystretch{1.2} 
 \begin{array}{|c|c|c|c|c|c|} 
 \hline 
 \mbox{Operad} & \mbox{Type of algebras} & \mbox{Koszul} &
  \mbox {Cyclic} & \mbox {Dihedral } & \mbox {Hopf} \\ 
\hline 
\Ass = \mathcal{LL}_1 = \Gass1& 
\mbox {associative} & {\rm yes} & {\rm 
 yes} & {\rm yes}&  {\rm yes} \\ 
 \hline 
\Poiss   = \mathcal{LL}_0 & \mbox {Poisson} & {\rm  yes} & {\rm 
 yes} & {\rm yes}&  {\rm  yes} \\ 
 \hline 
\mathcal{LL}_q,\ q \not= 0,\infty 
&\mbox {\rm $LL_q$-algebras} & {\rm yes } & {\rm 
 yes} & {\rm yes}&  {\rm yes}\\ 
\hline 
 \mathcal{LL}_\infty &\mbox {\rm $LL_\infty$-algebras} & {\rm yes } & {\rm 
 yes} & {\rm yes}& {\rm no}\\ 
  \hline 
 \Vinb =\Gass2 &{\rm Vinberg } & {\rm yes} & {\rm 
 no} & {\rm no} &  {\rm no}\\ 
 \hline 
\pre-Lie = \Gass3 & \mbox {\rm pre-Lie } & {\rm yes} & {\rm 
 no} & {\rm no} &  {\rm no}\\ 
 \hline 
 \Gass4& \mbox {\rm $G_4$-associative} & {\rm no} & {\rm 
 yes} & {\rm yes}&  {\rm no}\\ 
 \hline 
 \Gass5& \mbox {\rm $G_5$-associative} & {\rm no} & {\rm 
 no} & {\rm yes}& {\rm no}\\ 
 \hline 
 \Gass6& \mbox {\rm Lie-admissible} & {\rm yes} & {\rm 
 yes} & {\rm yes}& {\rm no}\\ 
 \hline 
 \end{array} 
\] 
\caption{\label{tab}%
Koszulness, quadratic cyclicity, dihedrality and Hopfness of operads
with one type~(3) operation.
}
\end{figure}

In Definition~\ref{dih} below we single out a property of quadratic operads 
responsible for the existence of the dihedral
cohomology~\cite{gelfand-manin:book,loday:AM87} 
of associated algebras. As far as we
know, this property has never been considered before.
Let $\calP = \Gamma(E)/(R)$ be a quadratic operad. Let $\lambda \in
\Sigma_2$ be the generator and define a left $\Sigma_2$-action on
$E$ using the operadic right $\Sigma_2$-action by $\lambda e := e
\lambda$, for $e \in E$. It follows from the universal property of
free operads that this action extends to a left 
$\Sigma_2$-action on $\Gamma(E)$.

\begin{definition}
\label{dih}
We say that a quadratic operad $\calP = \Gamma(E)/(R)$ is {\em
dihedral\/} if the left $\Sigma_2$-action on $\Gamma(E)$ induces a
left $\Sigma_2$-action on $\calP$.  A quadratic operad is {\em
cyclic dihedral\/}, if it is both cyclic and dihedral and if these two
structures are compatible, by which we mean that
\[
(\lambda u)\sigma = \lambda (u \sigma),
\]
for each $u \in \calP(n)$, $\lambda \in \Sigma_2$, $\sigma \in
\Sigmap_n$ and $n \geq 1$.
In other words, the cyclic and dihedral actions make each piece
$\calP(n)$ of a cyclic dihedral operad a left $\Sigma_2$-
right $\Sigmap_n$-bimodule. 
\end{definition}

\begin{remark}
{\rm 
We emphasize that dihedrality is a property defined only for quadratic
operads. We do not know how to extend this definition for a general
operad. Observe that the left $\Sigma_2$-action on $\Gamma(E)$ induces
an action on $\calP$ as required in Definition~\ref{dih} if and only
if the space of relations $R \subset \Gamma(E)(3)$ is $\Sigma_2$-stable.

The operad $\Gamma(V_{2,2})$ considered in Remark~\ref{pisi_na_letisti} is
quadratic, cyclic and dihedral, but not cyclic dihedral, because the
left $\Sigma_2$-action on $V_{2,2}$ is clearly not compatible with the
right $\Sigmap_2$-action. On the other hand, each cyclic quadratic
operad which is dihedral is cyclic dihedral. 
\/}
\end{remark}

We leave as an exercise to prove that all quadratic 
operads generated by one
operation of type (1) or (2) are dihedral. Therefore 
again the only interesting
case to investigate is a type (3) operation. 
The dihedrality is then easily understood if we write the axioms in 
the polarized form as follows.
Let $E = \KK[\Sigma_2]$ and decompose
\begin{equation}
\label{decomp}
\Gamma(E)(3) = \Gamma_+(E)(3) \oplus \Gamma_-(E)(3),
\end{equation}
where $\Gamma_+(E)(3)$ is the $\Sigma_3$-subspace of 
$\Gamma(E)(3)$ generated by compositions $x(yz)$ and $[x,[y,z]]$,
and $\Gamma_-(E)(3)$ is the $\Sigma_3$-subspace of 
$\Gamma(E)(3)$ generated by compositions $x[y,z]$ and $[x,yz]$.

In the pictorial language of Section~\ref{sec:monoidal-structures},
$\Gamma_+(E)(3)$ is the $\Sigma_3$-invariant subspace generated by
compositions of the following two types
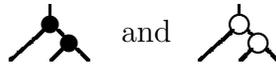
\[
\AAAA {}{}{} \hskip 1em \mbox { and } \hskip 1em \dvadva  {}{}{}
\]
while $\Gamma_-(E)(3)$ is the $\Sigma_3$-invariant subspace generated
by
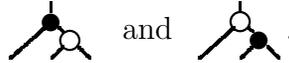
\[
\AAdva {}{}{} \hskip 1em \mbox { and } \hskip 1em \dvaAA  {}{}{}
\hskip .5em.
\]
Decomposition~(\ref{decomp}) is obviously $\Sigmap_3$-invariant. It is
almost evident that
$\lambda$ acts trivially on $\Gamma_+(E)(3)$ while on
$\Gamma_-(E)(3)$ it acts as the multiplication by $-1$. We therefore
get for free the following:

\begin{proposition}
\label{odmeny}
A quadratic operad $\calP = \Gamma(E)/(R)$ generated by a type~(3)
multiplication is dihedral if and only if the space of relations 
$R$ decomposes as
\[
R = R_+ \oplus R_- ,
\]
with  $R_+ \subset \Gamma_+(E)(3)$ and $R_- \subset \Gamma_-(E)(3)$.
\end{proposition}

\begin{example}
\label{nespravedliva_odmena}
{\rm
In this example we show that the operad $\Gass5$ for $G_5$-associative algebras
is dihedral but not cyclic. Recall from Example~\ref{Lie-adm}
that the polarized form of the
axioms for these algebras consists of the Jacobi identity
\[
 [x ,[ y , z]] + [y ,[ z , x]]+[z ,[ x , y]]=0
\]
and equation~(\ref{eq:7}) 
\[
 [x  y,z] +[y  z, x ]+[z x,y]=0. 
\]
Since the left-hand side of the Jacobi identity belongs to
$\Gamma_+(E)(3)$ and the right-hand side of~(\ref{eq:7})
to $\Gamma_-(E)(3)$, the space of relations obviously decomposes as
required by Proposition~\ref{odmeny}. Therefore $\Gass5$ is dihedral.

Let us inspect the cyclicity. By definition, the unpolarized form of
the axiom for $G_5$-associative algebras reads
\begin{equation}
\label{v_Mulhouse}
A(x,y,z) + A(y,z,x) + A(z,x,y) = 0,
\end{equation}
where $A$ denotes, as usual, the associator~(\ref{A}).
The action of $\gamma^+_3$ converts this equation to
\[
-A(x,y,z) + A(y,x,z) - A(y,z,x) = 0.
\]
The sum of the above equations gives 
\[
A(y,x,z) + A(z,x,y) = 0.
\]
It is then a simple linear algebra to prove that this equation does
not belong to the $\Sigma_3$-closure of~(\ref{v_Mulhouse}).
Therefore $\Gass5$ is not cyclic. 
}\end{example}

\begin{theorem}
\label{za_tyden_odjizdim}
Cyclic quadratic operads generated by one operation are dihedral.
\end{theorem}

{\bf Proof.}
The claim is obvious when $\calP$ is generated by one operation of
type (1) or (2). Suppose $\mathcal{P}$ is a quadratic operad of the form
\[
\mathcal{P}=\Gamma (E) /(R), \mbox { where } E =  :\KK[\Sigma_2].
\]
It was calculated in~\cite{getzler-kapranov:CPLNGT95} 
that, as $\Sigmap_3$-modules,
\begin{equation}
\label{jsem_utahany}
\Gamma_+(E)(3)=\id_3 \oplus V_{2,2} \oplus \sgn \oplus V_{2,2}
\ \mbox { and }\
\Gamma_-(E)(3)=V_{3,1} \oplus V_{2,1,1},
\end{equation}
where the irreducible representations $\id,\sgn, V_{2,2}, V_{3,1}$ and
$V_{2,1,1}$ are given by the following character table:

\[
\begin{array}{l|rrrrr}
& I & (01) & (012) & (0123) & (01)(23)\\
\hline 
\id & 1 & 1 & 1 & 1 & 1 \\
\sgn & 1 & -1 & 1 & -1 & 1 \\
V_{2,2} & 2 & 0 & -1 & 0 & 2 \\
V_{3,1} & 3 & 1 & 0 & -1 & -1 \\
V_{2,1,1} & 3 & -1 & 0 & 1 & -1 \\
\end{array}
\]

Observe that there are no common factors in
$\Gamma_+(E)(3)$ and
$\Gamma_-(E)(3)$, therefore  
it follows from an elementary representation theory
that each $\Sigmap_3$-invariant subspace $R$ of $\Gamma(E)(3)$
decomposes as $R = R_+ \oplus R_-$ 
with  $R_+ \subset \Gamma_+(E)(3)$ and $R_- \subset
\Gamma_-(E)(3)$. This means that $\calP$ is cyclic, by
Proposition~\ref{odmeny}.%
\qed

Theorem~\ref{za_tyden_odjizdim} was a consequence of the fact that for
operads generated by one operation, the $\Sigmap_3$-spaces
$\Gamma_+(E)(3)$ and $\Gamma_-(E)(3)$ do not contain a common
irreducible factor. The following example shows that this is not
longer true for general quadratic operads.

\begin{example}{\rm
Consider the quadratic operad
$\mathcal{P}=\Gamma(E)/(R)$, where $E := \KK[\Sigma_2] \oplus \id_2$
and where $(R)$ is 
the operadic ideal generated by the relations
\begin{eqnarray*}
r_1 & := &x \cdot (yz)+y \cdot (zx)+ z \cdot (xy)=0\
\mbox { and }
\\
r_2 & : =  &(xy) \cdot z + (z \cdot x) y +x (z \cdot y)=0.
\end{eqnarray*}
In the above display, $\cdot$ denotes the multiplication corresponding
to a generator of $\KK[\Sigma_2]$ and we, as usual, omit the symbol
for the commutative multiplication corresponding to a generator of $\id_2$.
Then $\mathcal{P}$ is cyclic but not dihedral.

Let us explain how this example was constructed.
It can be calculated that in
decomposition~(\ref{decomp}) of the 27-dimensional space
$\Gamma(E)(3)$,
\begin{eqnarray*}
\Gamma_+(E)(3) &=&  3 \id_2 \oplus \sgn \oplus 4 V_{2,2} \oplus V_{3,1}
\mbox { and } 
\\
\Gamma_-(E)(3) &=& 2V_{3,1} \oplus 2V_{2,1,1}.
\end{eqnarray*}
There is a common irreducible factor $V_{3,1}$ which occurs both in
$\Gamma_+(E)(3)$ and in $\Gamma_-(E)(3)$. Therefore, to construct 
an operad which is cyclic but not dihedral, it is
enough to choose a generator $e_+$ of $V_{3,1}$ in $\Gamma_+(E)$ and
a generator $e_-$ of $V_{3,1}$ in $\Gamma_-(E)$ and define $R$ to be
the $\Sigmap_3$-subspace of $\Gamma(E)(3)$ generated by $e_+ + e_-$. 
Operad $\calP$ above corresponds to one of these choices.
}\end{example}

In Figure~\ref{tab} we also recalled the following more or less well-known
results about Koszulness of operads considered in this paper.
 The operad $\Ass$ is Koszul by~\cite{ginzburg-kapranov:DMJ94} 
and the operad $\Poiss$ by~\cite[Corollary~4.6]{markl:dl}. The 
 operad $\Vinb$ is Koszul, because it is isomorphic to 
 the operad $\pre-Lie$ which is known to be
 Koszul~\cite{chapoton-livernet:pre-lie}. The operads $\Gass4$ and 
 $\Gass5$ are not Koszul, as proved in~\cite{remm:CRAS}. 

The Koszulness of the operad $\Gass6$ can be proved as follows. In
Example~\ref{Lie-adm} we
observed that $G_6$-associative algebras consist of a commutative
multiplication and a Lie bracket, with no relation between these two
operations. Therefore $\Gass6$ is the free product
\[
\Gass6 \cong \Lie * \Gamma(\id_2)
\]
of the operad $\Lie$ for Lie algebras and the free operad
$\Gamma(\id_2)$ generated by one commutative operation. The Koszulity
of the operad $\Gass6$ now follows from the obvious fact
that the free product of two quadratic Koszul operads is again
quadratic Koszul.

Let us turn our attention to the operad $\mathcal{LL}_q$ governing
$LL_q$-algebras. For $q = 0$, the Koszulness of 
$\mathcal{LL}_q$ follows from the isomorphism $\mathcal{LL}_0 \cong
\Poiss$. For $q \not\in \{0,\infty\}$ we argue as follows. The
Koszulness of an operad $\calP$ is characterized by 
the acyclicity, in positive dimensions, of the cobar dual of 
$\calP$~\cite{ginzburg-kapranov:DMJ94}. 
This means that Koszulness is not affected by a
field extension. We may therefore assume that the ground field $\KK$ is
algebraically closed. In this case the operad $\mathcal{LL}_q$ is
isomorphic, via the isomorphism~(\ref{za_chvili_prednaska}), 
to the operad $\Ass$ which is~\cite{markl:dl}. 

It remains to analyze the case $q = \infty$. It immediately follows
from the definition of $LL_\infty$-algebras that the corresponding
operad $\mathcal{LL}_\infty$ is constructed from Koszul quadratic operads
$\Gamma(\id_2)$ and $\Gamma(\sgn_2)/(\Gamma(\sgn_2)(3))$ via a
distributive law, it is therefore Koszul by~\cite[Theorem~4.5]{markl:dl}.


\def\cprime{$'$}

\catcode`\@=11

\noindent
Mathematical Institute of the Academy, 
\v Zitn\'a 25,
115 67 Prague 1,
The Czech Republic\hfill\break
e-mail: {\tt markl@math.cas.cz}

\noindent 
and

\noindent 
Laboratoire de Math\'ematiques et Applications,
Universit\'e de Haute Alsace, \hfill\break 
Facult\'e des Sciences et Techniques,
4, rue des Fr\`eres Lumi\`ere,
68093 Mulhouse cedex,
France \hfill\break 
e-mail: {\tt Elisabeth.Remm@uha.fr}

\vfill
\hfill 
{\tt \jobname.tex}

\end{document}